\documentclass{article}
\usepackage{amssymb,amsmath}
\usepackage{graphics}

\def\qed{\hfill {\hbox{${\vcenter{\vbox{               
   \hrule height 0.4pt\hbox{\vrule width 0.4pt height 6pt
   \kern5pt\vrule width 0.4pt}\hrule height 0.4pt}}}$}}}

\def\tr{\triangleright}

\newtheorem{theorem}{Theorem}
\newtheorem{definition}{Definition}
\newtheorem{lemma}[theorem]{Lemma}
\newtheorem{proposition}[theorem]{Proposition}
\newtheorem{corollary}[theorem]{Corollary}
\newtheorem{example}{Example}

\newtheorem{conjecture}{Conjecture}

\newenvironment{proof}[1][Proof]{\smallskip\noindent{\bf #1.}\quad}%
{\qed\par\medskip}

\author{{\begin{tabular}{c} Esteban Adam Navas \\ 
\small Department of Mathematics \\ 
\small University of California, Riverside\\
\small 900 University Avenue \\
\small Riverside, CA, 92521 \\
\small{\texttt{enava004@student.ucr.edu}}
\end{tabular}}
\and
{\begin{tabular}{c} Sam Nelson \\ \small
Department of Mathematics\\ 
\small Pomona College \\
\small 610 North College Avenue \\
\small Claremont, CA 91711 \\
\small{\texttt{knots@esotericka.org}}
\end{tabular}}
}

\date{}

\title{\Large \textbf{On symplectic quandles}}

\begin{document}
\maketitle

\begin{abstract}
We study the structure of \textit{symplectic quandles}, quandles which
are also $R$-modules equipped with an antisymmetric bilinear form. We show 
that every finite dimensional symplectic quandle over a finite field 
$\mathbb{F}$ or arbitrary field $\mathbb{F}$ of characteristic other than 2 
is 
a disjoint union of a trivial quandle and a connected quandle. We use the 
module structure of a symplectic quandle over a finite ring to refine and
strengthen the quandle counting invariant.
\end{abstract}



\section{\large \textbf{Introduction}}

A \textit{quandle} is a non-associative algebraic structure whose axioms
may be understood as transcriptions of the Reidemeister moves.
The term ``quandle'' was introduced by Joyce \cite{J}, though quandles 
have been studied by other authors under various names such as
 ``distributive groupoids'' \cite{M} and (for a certain special case) 
``Kei'' (\cite{T}, \cite{P}). Several 
generalizations of quandles have been defined and studied, including 
\textit{automorphic sets} (see \cite{B}) and \textit{racks} (see \cite{FR}) 
where the axioms are derived from regular 
isotopy moves, \textit{virtual quandles} (see \cite{V}) where additional 
structure is included for modeling virtual Reidemeister moves, and 
\textit{biquandles} and \textit{Yang-Baxter Sets}, which also have axioms 
derived from the Reidemeister moves but use a different correspondence 
between algebra elements and portions of link diagrams.

Quandles have found applications in topology as a source of invariants of 
topological spaces. In particular, finite quandles are useful for defining 
computable invariants of knotted circles in $S^3$ and other 3-manifolds 
as well as generalizations of ordinary knots such as virtual knots, knotted 
surfaces in $S^4$, etc.

In \cite{Y}, an example of a quandle structure defined on a module $M$ over a 
commutative ring $R$ with a choice of antisymmetric bilinear 
form $\langle,\rangle:M\times M\to R$ is given. In this paper we study the 
structure of this type of quandle, which we call a \textit{symplectic 
quandle}\footnote{After the completion of this paper, we were reminded that
symplectic quandles are also called \textit{quandles of transvections}.}. 
Our main result says that every  symplectic quandle $Q$ over a field 
$\mathbb{F}$ (of characteristic other than 2 if $\mathbb{F}$ is not finite) 
is \textit{almost connected}, that is, $Q$ is a disjoint 
union in the sense of \cite{B} of a trivial quandle and a connected quandle. 
Symplectic quandles are not just quandles but also $R$-modules; we show how 
to use the $R$-module structure of a finite symplectic quandle to 
enhance the usual quandle counting invariant.

The paper is organized as follows. In section \ref{qb} we recall the basic
definitions and standard examples of quandles. In section \ref{sq} we define
symplectic quandles, give some examples and show that symplectic quandles are 
almost connected. In section \ref{k} we give an application of symplectic 
quandles to knot invariants, defining a new family of enhanced
quandle counting invariants associated to finite symplectic quandles.

\section{\large \textbf{Quandle basics}}\label{qb}

We begin with a definition from \cite{J}.

\begin{definition}
\textup{Let $Q$ be a set and $\tr:Q\times Q\to Q$ a binary operation 
satisfying
\begin{list}{}{}
\item[(i)]{for all $a\in Q$, $a\tr a=a$,}
\item[(ii)]{for all $a,b\in Q$, there is a unique $c\in Q$ such that 
$a=c\tr b$, and}
\item[(iii)]{for all $a,b,c\in Q$, $(a\tr b)\tr c = (a\tr c)\tr(b\tr c)$.}
\end{list}}
\end{definition}

Axiom (ii) says that the quandle operation $\tr$ has a right inverse 
$\tr^{-1}$ such that $(x \tr y)\tr^{-1} y=x$ and $(x \tr^{-1} y)\tr y=x$. 
It is not hard to show that $Q$ is a quandle under $\tr^{-1}$ (called the 
\textit{dual} of $(Q,\tr)$) and that the two operations distribute over each 
other.

Standard examples of quandle structures include:

\begin{example}
\textup{Any set $Q$ is a quandle under the operation $x\tr y= x$, called a 
\textit{trivial} quandle. We denote the trivial quandle of order $n$ by $T_n$.}
\end{example}

\begin{example}
\textup{The finite abelian group $\mathbb{Z}_n$ is a quandle under 
$x\tr y=2y-x$. This is sometimes called the \textit{cyclic quandle} of order 
n.}
\end{example}

\begin{example}
\textup{Any group $G$ is a quandle under the following operations:
\begin{list}{$\bullet$}{}
\item{$x\tr y=y^{-1}xy$, or}
\item{$x\tr y=y^{-n}xy^n$, or}
\item{$x\tr y=s(xy^{-1})y$ where $s\in \mathrm{Aut}(G)$.}
\end{list}}
\end{example}

\begin{example}
\textup{Any module over $\mathbb{Z}[t^{\pm 1}]$ is a quandle under
\[x\tr y= tx+(1-t)y.\]
Quandles of this type are called \textit{Alexander quandles.} 
See \cite{AG} and
\cite{N} for more.}
\end{example}

\begin{example}
\textup{For any tame link diagram $L$, there is a quandle $Q(L)$ defined by a 
Wirtinger-style presentation with one generator for each arc and one relation
at each crossing. 
\[\includegraphics{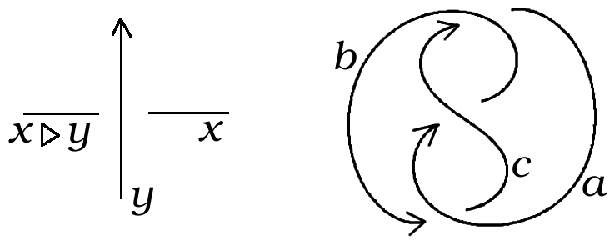} \ \ 
\raisebox{0.5in}{$Q(L)= \langle a,b,c \ | 
\ a\tr b=c, b\tr c=a, c\tr a =b \rangle.$}\]
This \textit{knot quandle} is in fact a classifying invariant
of knots and unsplit links in $S^3$ and certain other 3-manifolds up to
orientation-reversing homeomorphism of the ambient space. Elements of
a knot quandle are equivalence classes of quandle words in the arc generators 
under the equivalence relation generated by
the quandle axioms. See \cite{J} and \cite{FR} for more.}
\end{example}

\begin{definition}
\textup{Let $Q=\{x_1, x_2, \dots, x_n\}$ be a finite quandle. The matrix $M_Q$
with $M_Q[i,j]=k$ where $x_k=x_i\tr x_j$ for all $i,j\in \{1, 2,\dots, n\}$
is the \textit{quandle matrix} of $Q$. That is, $M_Q$ is the operation table of $Q$
without the ``$x$''s.}
\end{definition}

\begin{example}
\textup{The quandle $Q=\mathbb{Z}_3=\{1, 2, 3\}$ (note that we use 3 for 
the representative of the coset $0+3\mathbb{Z}$ so that our row and column 
numbers start with 1 instead of 0) with $i\tr j = 2j-i$ has quandle matrix
\[ M_Q=\left[\begin{array}{ccc}
1 & 3 & 2 \\
3 & 2 & 1 \\
2 & 1 & 3
\end{array}\right].\]}
\end{example}

\section{\large \textbf{Symplectic quandles}} \label{sq}

We begin this section with a definition (see \cite{Y}).

\begin{definition}\textup{
Let $M$ be a finite dimensional free module over a 
commutative ring with identity
$R$ and let $\langle,\rangle:M\times M\to R$ be an antisymmetric
bilinear form such that $\langle \mathbf{x},\mathbf{x}\rangle =0$ for all 
$\mathbf{x}\in M$. Then $M$ is a quandle with quandle operation}
\[ \mathbf{x}\tr \mathbf{y} = \mathbf{x}+\langle \mathbf{x},\mathbf{y}\rangle 
\mathbf{y}.\]
\textup{The dual quandle operation is given by }
\[ \mathbf{x}\tr^{-1} \mathbf{y} = 
\mathbf{x}-\langle \mathbf{x},\mathbf{y}\rangle \mathbf{y}.\]
\textup{If $R$ is a field and the form is non-degenerate, i.e., if
$\langle \mathbf{x},\mathbf{y} \rangle =\mathbf{0}$ for all
$\mathbf{y}\in M$ implies
$\mathbf{x}=\mathbf{0}\in M$, 
then $M$ is \textit{symplectic vector space}
and $\langle,\rangle$ is a \textit{symplectic form}; thus it is natural to 
refer to such $M$ as \textit{symplectic quandles}. For simplicity, we will 
use the term ``symplectic quandle over $R$'' to refer to the general case 
where $R$ is any ring and $\langle,\rangle$ is any antisymmetric bilinear 
form. 
If $\langle,\rangle$ is non-degenerate, we will say $(M,\tr)$ is a 
\textit{non-degenerate symplectic quandle over $R$}. $M$ and $M'$ are 
\textit{isometric} if there is an $R$-module isomorphism $\phi:M\to M'$ which 
preserves the bilinear form $\langle,\rangle$.}
\end{definition}

\begin{definition}\textup{
A quandle is \textit{involutory} if $\tr=\tr^{-1}$. Note that involutory 
quandles are also called \textit{kei} (see \cite{J}, \cite{P}  and \cite{T} 
for more).}
\end{definition}

\begin{proposition}
If $M$ is a symplectic quandle over a ring $R$ of characteristic 2, then
$M$ is involutory.
\end{proposition}

\begin{proof}
If $M$ is a symplectic quandle over a ring $R$ of characteristic 2, then
for any $\mathbf{x},\mathbf{y}\in M$ we have
\[ \mathbf{x}\tr \mathbf{y} 
= \mathbf{x}+\langle\mathbf{x},\mathbf{y}\rangle \mathbf{y} 
= \mathbf{x}-\langle\mathbf{x},\mathbf{y}\rangle \mathbf{y} 
=\mathbf{x}\tr^{-1} \mathbf{y}.\]
\end{proof}

\begin{definition}\label{du}
\textup{Let $Q$ and $Q'$ be quandles with $Q\cap Q'=\emptyset$. Then 
we can make $Q\cup Q'$ a quandle by defining $x\tr y = x$ when $x\in Q$
and $y\in Q'$ or when $x\in Q'$ and $y\in Q$. This is the \textit{disjoint union}
of $Q$ and $Q'$ in the sense of Brieskorn \cite{B}. If $Q$ and $Q'$ are 
finite then the matrix of $Q\cup Q'$ is the $(n+m) \times (n+m)$ block matrix
\[ M_{Q\cup Q'}=\left[\begin{array}{c|c}
M_Q & row \\ \hline
row & M_{Q'}
\end{array}\right]\]
where $row$ indicates that all entries are equal to their row number
and we denote $Q=\{x_1,\dots, x_n\},$ and $Q'=\{x_{n+1},\dots, x_{n+m}\}$.}
\end{definition}

Every quandle $Q$ can be decomposed as a disjoint union of a trivial 
subquandle \[D=\{ x\in Q \ | \ x\tr y =x \ \mathrm{and} \ y\tr x = y \ 
\forall y\in Q\}\] and a non-trivial subquandle $Q\setminus D$. Both $D$ 
and $Q\setminus D$ may be 
empty, and $Q\setminus D$ may contain trivial subquandles. Call $D$ the 
\textit{maximal trivial component} of $Q$.

\begin{example}
\textup{The quandle $Q$ with matrix $M_{Q}$ below has maximal trivial component
$D=\{x_5\}$ and $Q\setminus D=\{x_1, x_2, x_3, x_4\}$.}
\[M_{Q}=\left[\begin{array}{ccccc} 
1 & 1 & 1 & 2 & 1 \\
2 & 2 & 2 & 3 & 2 \\
3 & 3 & 3 & 1 & 3 \\
4 & 4 & 4 & 4 & 4 \\
5 & 5 & 5 & 5 & 5
\end{array}\right]\]
\textup{Notice that even though $\{x_1, x_2, x_3\}$ is a trivial subquandle,
it is not part of the maximal trivial component because of the way in which 
it is embedded in the overall quandle.}
\end{example}

\begin{proposition}
Let $Q$ be a symplectic quandle over $R$. Then the maximal trivial component of
$Q$ is the submodule of $R$ on which $\langle,\rangle$ is degenerate.
\end{proposition}

\begin{proof}
For any $\mathbf{x}\in Q$ we have
$\mathbf{x}\tr \mathbf{0} = 
\mathbf{x}+\langle\mathbf{x},\mathbf{0}\rangle\mathbf{0} =\mathbf{x} $
and
$\mathbf{0}\tr \mathbf{x} = 
\mathbf{0}+\langle\mathbf{0},\mathbf{x}\rangle\mathbf{x} =\mathbf{0},$
so $\mathbf{0}$ is in the maximal trivial component of $Q$. More generally, 
let $D$ be the submodule of $Q$ on which $\langle,\rangle$ is 
degenerate, i.e., 
\[D=\{ \mathbf{x}\in Q \ | \ \langle \mathbf{x}, \mathbf{y}\rangle =0 
\ \forall \mathbf{y}\in Q\}.\]
Then for any $\mathbf{d}\in D$ we have $\mathbf{x}\tr \mathbf{d} = 
\mathbf{x}+0\mathbf{d} = \mathbf{x}$ and $\mathbf{d}\tr \mathbf{x} = 
\mathbf{d}+ 0\mathbf{x} = \mathbf{d}$,
so $D$ is a trivial subquandle of $Q$ and $Q$ is the disjoint union of $D$ and 
$Q\setminus D$ in sense of definition \ref{du}. If $\mathbf{x}\not\in D$, 
then there 
is some $\mathbf{y}\in Q$ with $\langle\mathbf{x},\mathbf{y}\rangle\ne0$ 
so that $\mathbf{x}\tr \mathbf{y} \ne \mathbf{x}$;
then $Q\setminus D$ is non-trivial and $D$ is precisely the submodule of $Q$
on which $\langle,\rangle$ is degenerate.
\end{proof}

\begin{corollary}
If $Q$ is a nondegenerate symplectic quandle, then the maximal trivial component 
of $Q$ is $D=\{\mathbf{0}\}$.
\end{corollary}

We will now restrict our attention to the case where $M$ is a free module
over a PID $R$. It is a standard result (see \cite{AW} for example) that such 
an $M$ equipped with a nondegenerate antisymmetric bilinear form must be even 
dimensional, with basis
$\{\mathbf{b}_i \ | \ i=1,\dots, 2n\}$ such that
\[\left\langle \mathbf{x},\mathbf{y}\right\rangle =
\left\langle \sum_{i=1}^{2n}x_i\mathbf{b}_i,
\sum_{i=1}^{2n}y_i\mathbf{b}_i\right\rangle =
\sum_{i=1}^{2n} \epsilon(i) \alpha_i x_i y_{i+\epsilon(i)}
\quad 
\mathrm{where} \quad 
\epsilon(i)=\left\{\begin{array}{rl} 1 & i\ \mathrm{odd} \\
-1 & i \ \mathrm{even},\end{array}\right. \]
$\alpha_{2i}=\alpha_{2i-1}$,  
and each $\alpha_i$ is either $1$ or a nonunit in $R$. Such a basis is called 
a \textit{symplectic basis}. The $\alpha_i$s 
are called \textit{invariant factors}, and the set with multiplicities of 
invariant factors determines the symplectic module structure up to isometry 
(i.e., $\langle,\rangle$-preserving isomorphism of $R$-modules). In particular, 
if $R$ is a field, then we may choose our basis so that 
$\alpha_{2i}=\alpha_{2i-1}=1$ for all $i=1,\dots, n$. 

In matrix notation with $\mathbf{x},$ $\mathbf{y}$ row vectors, we have 
$ \langle\mathbf{x},\mathbf{y}\rangle=\mathbf{x}A\mathbf{y}^T$
where $A$ is a block diagonal matrix of the form
\[A=\left[\begin{array}{cc|cc|c|cc}
0 & \alpha_2  & 0 & 0 & \dots & 0 & 0 \\ 
-\alpha_2 & 0 & 0 & 0 & \dots & 0 & 0 \\ \hline
0 & 0 & 0 & \alpha_4  & \dots & 0 & 0 \\ 
0 & 0 & -\alpha_4 & 0 & \dots & 0 & 0 \\ \hline
\vdots & \vdots & \vdots & \vdots & \ddots & \vdots & \vdots \\ \hline
0 & 0 & 0 & 0 & \dots & 0 & \alpha_{2n} \\
0 & 0 & 0 & 0 & \dots & -\alpha_{2n} & 0 
\end{array}\right].\]

It is clear that isometric $R$-modules are isomorphic as quandles. Conversely,
a symplectic quandle structure on $R^n$ determines the antisymmetric bilinear 
form $\langle,\rangle$ uniquely up to choice of basis: for a 
basis $\{\mathbf{b}_i \ | \ i=1,\dots, 2n\}$ of $R^{2n}$ we have
\[ \mathbf{b}_i \tr \mathbf{b}_j - \mathbf{b}_i=
\langle \mathbf{b}_i,\mathbf{b}_j\rangle \mathbf{b}_j 
=\alpha_{ij}\mathbf{b}_j\]
and since $\{\mathbf{b}_i\}$ is a basis, the $\alpha_{ij}$ thus determined is 
unique. Changing bases to get a symplectic basis, we then obtain the invariant
factors. Thus the quandle structure together with the $R$-module structure of $M$ 
determine the invariant factors and hence determine $\langle,\rangle$, and we 
have:

\begin{theorem}
Let $Q$ and $Q'$ be non-degenerate $2n$-dimensional symplectic quandles over a 
PID $R$. Then $Q$ and $Q'$ are isomorphic as quandles if and only if they are 
isometric.
\end{theorem}

Our search through examples of finite symplectic quandles of small cardinality
over $\mathbb{Z}_n$ for $n$ non-prime has failed to yield any examples of 
symplectic quandles which are isomorphic as quandles but not isometric as 
$R$-modules. Thus, we have

\begin{conjecture}\label{c:m}
Two symplectic quandles of the same dimension over $\mathbb{Z}_n$ are 
isomorphic as quandles if and only if they are isometric.
\end{conjecture}

The following example shows that cardinality alone does not determine $R$ or 
the rank of $Q$.

\begin{example}\label{e16}
\textup{Let $R=\mathbb{Z}_2$ and $F=\mathbb{Z}_2[t]/(t^2+t+1)$. Both $R$ and
$F$ are fields of characteristic 2, and the symplectic vector spaces
\[V=R^4, \quad \langle\mathbf{x},\mathbf{y}\rangle = 
\mathbf{x}\left[\begin{array}{cccc}
0 & 1 & 0 & 0 \\
1 & 0 & 0 & 0 \\
0 & 0 & 0 & 1 \\
0 & 0 & 1 & 0
\end{array}\right] \mathbf{y}^T\]
and 
\[V'=F^2, \quad \langle\mathbf{x},\mathbf{y}\rangle = 
\mathbf{x}\left[\begin{array}{cc}
0 & 1 \\
1 & 0 \\
\end{array}\right] \mathbf{y}^T\]
are both symplectic quandles of order 16. From their quandle
matrices
\[
M_{V}=\left[\begin{array}{cccccccccccccccc}
1  & 1  & 1  & 1  & 1  &  1 &  1 &  1 &  1 &  1 &  1 &  1 &  1 &  1 &  1 &  1 \\
2  & 2  & 4  & 3  & 2  &  2 &  8 &  7 &  2 &  2 & 12 & 11 &  2 &  2 & 16 & 15 \\
3  & 4  & 3  & 2  & 3  &  8 &  3 &  6 &  3 & 12 &  3 & 10 &  3 & 16 &  3 & 14 \\
4  & 3  & 2  & 4  & 4  &  7 &  6 &  4 &  4 & 11 & 10 &  4 &  4 & 15 & 14 &  4 \\
5  & 5  & 5  & 5  & 5  &  5 &  5 &  5 & 13 & 14 & 15 & 16 &  9 & 10 & 11 & 12 \\
6  & 6  & 8  & 7  & 6  &  6 &  4 &  3 & 14 & 13 &  6 &  6 & 10 &  9 &  6 &  6 \\
7  & 8  & 7  & 6  & 7  &  4 &  7 &  2 & 15 &  7 & 13 &  7 & 11 &  7 &  9 &  7 \\
8  & 7  & 6  & 8  & 8  &  3 &  2 &  8 & 16 &  8 &  8 & 13 & 12 &  8 &  8 &  9 \\
9  & 9  & 9  & 9  & 13 & 14 & 15 & 16 &  9 &  9 &  9 &  9 &  5 &  6 &  7 &  8 \\
10 & 10 & 12 & 11 & 14 & 13 & 10 & 10 & 10 & 10 &  4 &  3 &  6 &  5 & 10 & 10 \\
11 & 12 & 11 & 10 & 15 & 11 & 13 & 11 & 11 &  4 & 11 &  2 &  7 & 11 &  5 & 11 \\
12 & 11 & 10 & 12 & 16 & 12 & 12 & 13 & 12 &  3 &  2 & 12 &  8 & 12 & 12 &  5 \\
13 & 13 & 13 & 13 & 9  & 10 & 11 & 12 &  5 &  6 &  7 &  8 & 13 & 13 & 13 & 13 \\
14 & 14 & 16 & 15 & 10 &  9 & 14 & 14 &  6 &  5 & 14 & 14 & 14 & 14 &  4 &  3 \\
15 & 16 & 15 & 14 & 11 & 15 &  9 & 15 &  7 & 15 &  5 & 15 & 15 &  4 & 15 &  2 \\
16 & 15 & 14 & 16 & 12 & 16 & 16 & 9  &  8 & 16 & 16 &  5 & 16 &  3 &  2 & 16 
\end{array}\right]
\] and
\[
M_{V'}=\left[\begin{array}{cccccccccccccccc}
 1 &  1 &  1 &  1 &  1 &  1 &  1 &  1 &  1 &  1 &  1 &  1 &  1 &  1 &  1 &  1 \\
 2 &  2 &  2 &  2 &  6 & 5 &  8 &  7 & 14 & 16 & 15 & 13 & 10 & 11 &  9 & 12 \\
 3 &  3 &  3 &  3 & 11 &  9 & 10 & 12 &  7 &  6 &  8 &  5 & 15 & 16 & 13 & 14 \\
 4 &  4 &  4 &  4 & 16 & 13 & 15 & 14 & 12 & 11 & 10 &  9 &  8 &  6 &  5 &  7 \\
 5 &  6 &  8 &  7 &  5 &  2 & 16 & 11 &  5 & 14 & 12 &  3 &  5 & 10 &  4 & 15 \\
 6 &  5 &  7 &  8 &  2 &  6 &  9 & 13 & 10 &  3 &  6 & 15 & 14 &  4 & 12 &  6 \\
 7 &  8 &  6 &  5 & 15 & 10 &  7 &  2 &  3 &  9 & 13 &  7 & 11 &  7 & 16 &  4 \\
 8 &  7 &  5 &  6 & 12 & 14 &  2 &  8 & 16 &  8 &  3 & 11 &  4 & 13 &  8 &  9 \\
 9 & 11 & 10 & 12 &  9 &  3 &  6 & 16 &  9 &  7 & 14 &  4 &  9 & 15 &  2 &  8 \\
10 & 12 &  9 & 11 & 14 &  7 &  3 & 10 &  6 & 10 &  4 & 16 &  2 &  5 & 10 & 13 \\
11 &  9 & 12 & 10 &  3 & 11 & 13 &  5 & 15 &  4 & 11 &  8 &  7 &  2 & 14 & 11 \\
12 & 10 & 11 &  9 &  8 & 15 & 12 &  3 &  4 & 13 &  5 & 12 & 16 & 12 &  6 &  2 \\
13 & 16 & 15 & 14 & 13 &  4 & 11 &  6 & 13 & 12 &  7 &  2 & 13 &  8 &  3 & 10 \\
14 & 15 & 16 & 13 & 10 &  8 & 14 &  4 &  2 &  5 &  9 & 14 &  6 & 14 & 11 &  3 \\
15 & 14 & 13 & 16 &  7 & 12 &  4 & 15 & 11 & 15 &  2 &  6 &  3 &  9 & 15 &  5 \\
16 & 13 & 14 & 15 &  4 & 16 &  5 &  9 &  8 &  2 & 16 & 10 & 12 &  3 &  7 & 16
\end{array}\right]
\] we can easily see that $V$ and $V'$ are not isomorphic as quandles
by checking that the quandle polynomials $qp_{V}(s,t)=s^{16}t^{16}+15s^8t^8$ and
$qp_{V'}(s,t)=s^{16}t^{16}+15s^4t^4$ are not equal (see \cite{N2} for more). 
}
\end{example}

\begin{definition}
\textup{A quandle $Q$ is \textit{connected} if it has a single orbit, i.e.,
if every element $z\in Q$ can be obtained from every other element $x\in Q$ by
a sequence of quandle operations $\tr$ and dual quandle operations $\tr^{-1}$. 
A quandle is
\textit{almost connected} if it is a disjoint union in the sense of definition
\ref{du} of its maximal trivial component and a single connected subquandle.}
\end{definition}

Our main result says that symplectic quandles $Q$ over a finite field or
infinite field $\mathbb{F}$ of characteristic other than 2 are almost connected; 
in particular, if $\langle,\rangle$ is nondegenerate then the subquandle 
$Q\setminus \{\mathbf{0}\}$ is a connected 
quandle. Connected quandles are of particular interest for defining knot 
invariants since knot quandles for knots (i.e., single-component links) are always 
connected. In particular, the image of a quandle homomorphism 
$f:Q(L)\to T$ from a knot quandle to $T$ always lies within a single orbit 
of the codomain quandle $T$, though of course $f$ need not be surjective.

For the remainder of this section, let $Q$ be a symplectic quandle over a field 
$\mathbb{F}$ and choose a symplectic basis $\{\mathbf{b}_i\}$ with invariant 
factors $\alpha_{2i}=1$ for $i=1, \dots, n$.

\begin{lemma} \label{l:4}
If any component $x_i$ of 
$\displaystyle{\mathbf{x}=\sum_{i=1}^{2n} x_i\mathbf{b}_i \in Q}$ is 
nonzero then for any $j\in \{1,\dots, 2n\}$ there is a 
$\mathbf{z}=\mathbf{x}\tr \mathbf{y}\in Q$ with $z_j\ne 0$ for some $\mathbf{y}\in Q$. 
That is, we can change 
a zero component to a nonzero component using a quandle operation, provided at least 
one other component of $\mathbf{x}$ is nonzero.
\end{lemma}

\begin{proof}
Suppose $x_i\ne 0$ and $x_j=0$. Then choose $\beta \in \mathbb{F}$ such that
$\displaystyle{\beta \ne \frac{\epsilon(i)x_i}{\epsilon(j)x_{j+\epsilon(j)}}}$
or, if $x_{j+\epsilon(j)}=0,$ $\beta \ne -\epsilon(i)x_i,$
and define $\mathbf{y}= 
\mathbf{b}_{i+\epsilon(i)}+\beta\mathbf{b}_{j}$. Then we have
\begin{eqnarray*}
\mathbf{x}\tr \mathbf{y} & = & 
\mathbf{x}+ 
\left(\epsilon(i)x_i-\epsilon(j)x_{j+\epsilon(j)}\beta\right) \mathbf{y}
\end{eqnarray*}
and the $j$th component of $\mathbf{z}=\mathbf{x}\tr \mathbf{y}$ is 
\[z_j=0+(\epsilon(i)x_i-\beta\epsilon(j)x_{j+\epsilon(j)})\beta\]
which is nonzero by our choice of $\beta$.
\end{proof}

\begin{lemma} \label{l:5}
For any $\mathbf{x}=\displaystyle{\sum_{i=1}^{2n} x_i\mathbf{b}_i}\in Q$
and for any $\beta\in\mathbb{F}$, we 
can add (or subtract) $\beta^2x_i$ to (or from) $x_{i+\epsilon(i)}$ with 
quandle operations and dual quandle operations.
\end{lemma}

\begin{proof}
\[\mathbf{x} \tr \beta \mathbf{b}_{i+\epsilon(i)} 
= \mathbf{x}+ (\epsilon(i)x_{i}\beta)\beta \mathbf{b}_{i+\epsilon(i)} 
= \mathbf{x}+ \epsilon(i)\beta^2x_i \mathbf{b}_{i+\epsilon(i)} 
\]
and similarly
\[\mathbf{x} \tr^{-1} \beta \mathbf{b}_{i+\epsilon(i)} =
\mathbf{x}- \epsilon(i)\beta^2x_i \mathbf{b}_{i+\epsilon(i)}. \]
\end{proof}

\begin{lemma}\label{l:6} 
If the characteristic of $\mathbb{F}$ is not 2, then
for any $\mathbf{x}\ne 0$, we can change any component $x_i$ of $\mathbf{x}$ 
to any value $z\in\mathbb{F}$ with quandle operations and dual quandle operations.
\end{lemma}

\begin{proof}
Write $x_i=z+w$. By lemma \ref{l:4} we may assume that 
$x_{i+\epsilon(i)}\ne 0$.
Then 
\[\mathbf{x}\tr w\mathbf{b}_{i}
= \mathbf{x}+(\epsilon(i+\epsilon(i))x_{i+\epsilon(i)}w) w\mathbf{b}_{i} \]
and the new quandle element has $i$th component equal to
\begin{eqnarray*}
x_i+\epsilon(i+\epsilon(i))x_{i+\epsilon(i)}w^2 & = 
& z+w+\epsilon(i+\epsilon(i))x_{i+\epsilon(i)}w^2 \\
 & = & z+ \epsilon(i+\epsilon(i))x_{i+\epsilon(i)}\left( 
 \frac{w}{\epsilon(i+\epsilon(i))x_{i+\epsilon(i)}} + w^2\right).
\end{eqnarray*}
Let us denote $j=i+\epsilon(i)$. 
If the characteristic of $\mathbb{F}$ is not 2, then we can complete the square
to obtain
\begin{eqnarray*}
x_i+\epsilon(j)x_{j}w^2 & = & 
z+ \epsilon(j)x_{j}\left( \frac{1}{4x_{j}^2}+
\frac{w}{\epsilon(j)x_{j}} + w^2\right) 
-\epsilon(j)x_{j}\frac{1}{4x_{j}^2} \\ 
& = &
z+\epsilon(j) x_{j}\left(\frac{1}{2x_{j}}+\epsilon(j)w\right)^2 -
\epsilon(j) x_{j}\left(\frac{1}{2x_{j}}\right)^2.
\end{eqnarray*}
Then by lemma \ref{l:5} we can remove both terms via quandle operations
and dual quandle operations
to obtain $z$ in the $i$th component, as required.
\end{proof}

\begin{lemma}\label{l:7}
In a finite field $\mathbb{F}$ of characteristic 2, every element of 
$\mathbb{F}$ is a square.
\end{lemma}

\begin{proof}
The map $f:\mathbb{F}\to \mathbb{F}$ given by $f(x)=x^2$ is a homomorphism
of fields since \[f(x+y)=(x+y)^2=x^2+2xy+y^2=x^2+y^2=f(x)+f(y) \quad \mathrm{and}
\quad f(xy)=(xy)^2=x^2y^2.\]
Then $\mathrm{ker}(f)=\{0\}$ since $\mathbb{F}$ has no zero divisors; thus
$f$ is injective and, since $\mathbb{F}$ is finite, surjective. In particular, 
every $\alpha\in \mathbb{F}$ satisfies $\alpha=\beta^2$ for some 
$\beta\in \mathbb{F}$.
\end{proof}

Taken together, Lemmas \ref{l:4}, \ref{l:5}, \ref{l:6} and \ref{l:7} imply:

\begin{theorem}
Let $\mathbb{F}$ be a field of characteristic other than 2, or a finite
field of characteristic 2. Then every symplectic quandle over $\mathbb{F}$
is almost connected.
\end{theorem}

If $R$ is not a field, then symplectic quandles over $R$ need not be almost 
connected, as the next example shows.

\begin{example}
\textup{The symplectic quandle $V''=(\mathbb{Z}_4)^2$ with 
bilinear form}
\[\langle\mathbf{x},\mathbf{y}\rangle=
\mathbf{x}
\left[\begin{array}{cc}
0 & 2 \\
2 & 0
\end{array}\right] \mathbf{y}^T\]
\textup{has quandle matrix 
\[
M_{V''}=\left[\begin{array}{cccccccccccccccc}
 1 &  1 &  1 &  1 &  1 &  1 &  1 &  1 &  1 &  1 &  1 &  1 &  1 &  1 &  1 &  1 \\
 2 &  2 &  2 &  2 & 10 & 12 & 10 & 12 &  2 &  2 &  2 &  2 & 10 & 12 & 10 & 12 \\
 3 &  3 &  3 &  3 &  3 &  3 &  3 &  3 &  3 &  3 &  3 &  3 &  3 &  3 &  3 &  3 \\
 4 &  4 &  4 &  4 & 12 & 10 & 12 & 10 &  4 &  4 &  4 &  4 & 12 & 10 & 12 & 12 \\
 5 &  7 &  5 &  7 &  5 & 15 &  5 & 15 &  5 &  7 &  5 &  7 &  5 & 15 &  5 & 15 \\
 6 &  8 &  6 &  8 & 14 &  6 & 14 &  6 &  6 &  8 &  6 &  8 & 14 &  6 & 14 &  6 \\
 7 &  5 &  7 &  5 &  7 & 13 &  7 & 13 &  7 &  5 &  7 &  5 &  7 & 13 &  7 & 13 \\
 8 &  6 &  8 &  6 & 16 &  8 & 16 &  8 &  8 &  6 &  8 &  6 & 16 &  8 & 16 &  8 \\
 9 &  9 &  9 &  9 &  9 &  9 &  9 &  9 &  9 &  9 &  9 &  9 &  9 &  9 &  9 &  9 \\
10 & 10 & 10 & 10 &  2 &  4 &  2 &  4 & 10 & 10 & 10 & 10 &  2 &  4 &  2 &  4 \\
11 & 11 & 11 & 11 & 11 & 11 & 11 & 11 & 11 & 11 & 11 & 11 & 11 & 11 & 11 & 11 \\
12 & 12 & 12 & 12 &  4 &  2 &  4 &  2 & 12 & 12 & 12 & 12 &  4 &  2 &  4 &  2 \\
13 & 15 & 13 & 15 & 13 &  7 & 13 &  7 & 13 & 15 & 13 & 15 & 13 &  7 & 13 &  7 \\
14 & 16 & 14 & 16 &  6 & 14 &  6 & 14 & 14 & 16 & 14 & 16 &  6 & 14 &  6 & 14 \\
15 & 13 & 15 & 13 & 15 &  5 & 15 &  5 & 15 & 13 & 15 & 13 & 15 &  5 & 15 &  5 \\
16 & 14 & 16 & 14 &  8 & 16 &  8 & 16 & 16 & 14 & 16 & 14 &  8 & 16 &  8 & 16  
\end{array}\right].
\]
$V''$ has maximal trivial component $D=\{x_1,x_3,x_9,x_{11}\}$, but the nontrivial 
component $V''\setminus D$ has 
disjoint orbit subquandles $\{x_2,x_4,x_{10},x_{12}\}$, $\{x_5,x_7,x_{13},x_{15}\}$ 
and $\{x_6,x_8,x_{14},x_{16}\}$ 
and hence is not connected. For comparison with the order 16 symplectic 
quandles in example \ref{e16}, the quandle polynomial for $V''$
is $qp_{V''}(s,t)=4s^{16}t^{16}+12s^8t^8.$}
\end{example}

\section{\large \textbf{Symplectic quandles and knot invariants}} \label{k}

The primary application for finite quandles has so far been in the construction
of link invariants. Given a finite quandle $T$ we have the quandle counting
invariant $|\mathrm{Hom}(Q(L),T)|$, the quandle 2-cocycle invariants 
$\Phi_{\chi}(L,T)$ and the specialized subquandle polynomial invariants
$\Phi_{qp}(L)$ described in \cite{C4} and \cite{N2} respectively. The 
connected component of a symplectic quandle over a finite field is a 
finite connected quandle which generally has a number of nontrivial 
subquandles, making this type of quandle well suited for the specialized 
subquandle polynomial invariant. In this section we describe two additional
ways of getting extra information about the knot or link type from the set 
of homomorphisms from a link quandle into a finite symplectic quandle.

One easy way to get more information out of the set $|\mathrm{Hom}(Q(L),T)|$
is to count the cardinalities of the image subquandles for each
$f\in \mathrm{Hom}(Q(L),T)$; even if $T$ is connected, the smallest subquandle
of $T$ containing the images of generators of $Q(L)$ need not be the entire
quandle $T$. If instead of counting 1 for each homomorphism $f$, we 
count the cardinality of the image of $f$, we obtain a set with multiplicities 
of integers, which we can convert into a polynomial for easy comparison with 
other invariant values by converting the elements of the set to
exponents of a variable $q$ and converting the multiplicities to coefficients.
Thus we have

\begin{definition}
\textup{The \textit{enhanced quandle counting invariant} of a link $L$ with
respect to a finite target quandle $T$ is given by}
\[\Phi_E(L,T)=\sum_{f\in \mathrm{Hom}(Q(L),T)} q^{|\mathrm{Im}(f)|}.\]
\end{definition}

This enhanced quandle counting invariant can be understood as a decomposition
of the usual quandle counting invariant into a sum of counting invariants
over all subquandles of our target quandle with the restriction that we only
count \textit{surjective} homomorphisms onto each subquandle.

For any subquandle $S\subseteq T$ of a finite quandle $T$,
let $\mathrm{SH}(Q(L),S)$ be the set of surjective quandle homomorphisms
from a link quandle $Q(L)$ onto $S$ and let $SQ(T)$ be the set of all subquandles
of $T$. Then
\[\Phi_E(L,T) = \sum_{S\in SQ(T)} |\mathrm{SH}(Q(L),S)|q^{|S|}.\]

Because symplectic quandles are not just quandles but also $R$-modules, we
can take advantage of the $R$-module structure of a finite symplectic quandle
$T$ to 
further enhance the counting invariant.

\begin{definition}\label{sqp}\textup{
Let $T$ be a finite symplectic quandle over a (necessarily finite) ring $R$ 
and let $L$ be a link. Then for each $f\in \mathrm{Hom}(Q(L),T)$, let 
$\rho(f)$ be the cardinality of the $R$-submodule spanned by 
$\mathrm{Im}(f)\subseteq T$ (note that $\mathrm{Im}(f)$ itself 
need not be a submodule).
Then the \textit{symplectic quandle polynomial} of $L$ with respect to $T$ is}
\[\Phi_{sqp}(L,T)=\sum_{f\in \mathrm{Hom}(Q(L),T)} 
q^{|\mathrm{Im}(f)|}z^{\rho(f)}.\]
\end{definition}

Note that in definition \ref{sqp} the finite target quandle $T$ has a fixed
$R$-module structure; in the case of a counterexample to conjecture \ref{c:m},
i.e., if two symplectic quandles exist which are isomorphic as quandles but 
not as modules, then we would expect two such symplectic quandles to define 
distinct symplectic quandle polynomial invariants. In particular, if $R$ is
not a field then we must be careful to specify the $R$-module
structure of $T$ and our choice of bilinear form.

The following example demonstrates that $\Phi_{sqp}$ contains more information 
than the quandle counting invariant alone.

\begin{example}
\textup{The two pictured virtual links have the same value for the
quandle counting invariant with respect to the symplectic quandle 
$T=(\mathbb{Z}_3)^2$ but different values for $\Phi_{sqp}(L,T)$. }
\[\begin{array}{c}
\raisebox{0.5in}{$L_1:$}\includegraphics{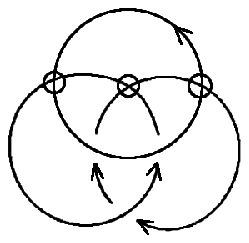} \\
|\mathrm{Hom}(Q(L_1),T)|=105 \\
\Phi_{sqp}(L_1,T)=9qz+72q^2z^3+24qz^3 \end{array} \quad \quad
\begin{array}{c}
\raisebox{0.5in}{$L_2:$}\includegraphics{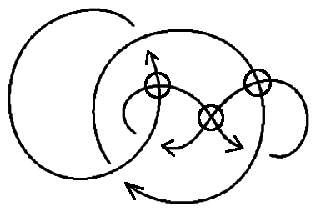} \\
|\mathrm{Hom}(Q(L_2),T)|=105 \\ 
\Phi_{sqp}(L_2,T)=qz+72q^2z^3+24q^3z^3+8qz^3\end{array}
\]
\end{example}

\begin{proposition}
If $T=K^{2n}$ is the nondegenerate symplectic quandle of 
dimension $2n$ over the Galois field $K=GF(p^m)$ for a prime $p$, then 
$\Phi_{sqp}(\mathrm{Unknot},T)=qz+(p^{2nm}-1)qz^{p^m}$.
\end{proposition}

\begin{proof}
Every element of $\mathrm{Hom}(Q(\mathrm{Unknot}),T)$ is a constant
map into a single element of $T$. The zero map contributes $q^1z^1=qz$
to the sum, while each of the nonzero constant maps has image subquandle 
consisting of a single element of $T$ which spans a dimension 1 subspace;
hence each of these $p^{2nm}-1$ maps contributes $qz^{p^m}$ to the sum. 
\end{proof}

Specializing $z=1$ and $q=1$ in $\Phi_{sqp}(L,T)$ yields the quandle counting
invariant $|\mathrm{Hom}(Q(L),T)|.$ Specializing $z=1$ yields the enhanced 
quandle counting invariant $\Phi_E(L,T)$.

Our initial computations suggest that these symplectic quandle polynomial 
invariants are quite non-trivial for virtual links, though the fact that finite 
symplectic quandles tend to have rather large cardinality ($|R|^{2n}$) means 
that more efficient computing algorithms may be required to explore these 
invariants in greater detail. Our \texttt{Maple} software is able to compute 
$\Phi_{sqp}(L,T)$ for links with smallish numbers of crossings for symplectic 
quandles of order $\le 81$ in a relatively short amount of time, but the time 
requirement increases rapidly as $|R|$ and $n$ increase.

\end{document}